\input ./style/arxiv-general.cfg
\documentclass[MSNbibl,number,citesort,dvips]{arxbj}
\makeatletter
   \@ifpackageloaded{graphicx}{}{\usepackage{graphicx}}
\makeatother
\aid{0}
\volume{22}
\issue{1}
\pubyear{2016}
\firstpage{563}
\lastpage{588}
\doi{10.3150/14-BEJ668} % kopijuoti is laisko
\docsubty{FLA}

\makeatletter

\newtheorem{them}{Theorem}
\newtheorem{lemma}{Lemma}%[section]
\newtheorem{prop}[them]{Proposition}
\newtheorem{cor}[them]{Corollary}
\newremark{Remark}{Remark}
\newcommand{\Rd}{\mathbb{R}^d}
\newcommand{\eqref}[1]{(\ref{#1})}
\renewcommand{\epsilon}{\varepsilon}

\makeatother

\begin{document}
\begin{frontmatter}

\title{Unitary transformations, empirical processes and distribution free testing}
\runtitle{Unitary transformations of empirical processes}

\begin{aug}
%%%% inicialai - be tarpu
\author{\inits{E.}\fnms{Estate} \snm{Khmaladze}\ead
[label=e1]{Estate.Khmaladze@vuw.ac.nz}}% \and
%\author{\inits{}\fnms{}~\snm{}\thanksref{}\ead[label=e2]{}}
%\author{\inits{}\fnms{}~\snm{}}
%%\runauthor{} %% auto
%\dedicated{}
\address[]{School of Mathematics, Statistics and Operations Research, Victoria University of Wellington, Wellington, New Zealand. \printead{e1}}
%\address[]{}
\end{aug}

% HISTORY:
\received{\smonth{11} \syear{2013}}
\revised{\smonth{7} \syear{2014}}

% ABSTRACT
%
\begin{abstract}
The main message in this paper is that there are surprisingly many
different Brownian bridges, some of them -- familiar, some of them --
less familiar. Many of these Brownian bridges are very close to
Brownian motions. Somewhat loosely speaking, we show that all the
bridges can be conveniently mapped onto each other, and hence, to one
``standard'' bridge.

The paper shows that, a consequence of this, we obtain a unified theory
of distribution free testing in $\Rd$, both for discrete and
continuous cases, and for simple and parametric hypothesis.
\end{abstract}

% KEYWORDS
% visi is mazosios raides ir pagal abecele
%
\end{frontmatter}

%s1 #&#
\section{Introduction}\label{sec1}
Let $X_1,X_2,\dots,X_n$ be a sequence of i.i.d. random variables in
$\Rd$ with distribution $F$, and consider an empirical process based
on this sequence:
\[
v_{nF}(B)=\sqrt{n} \bigl[F_n(B) - F(B)\bigr],
\]
where $B$ is a Borel subset of $\Rd$ and
\[
F_n(B)=\frac{1}{n} \sum_{i=1}^n
I_{\{X_i\in B\}}
\]
is an empirical distribution. If sets $B$ are chosen as unbounded
rectangles $(-\infty, x] = (-\infty, x_1]\times\cdots\times(-\infty
, x_d]$, then we obtain more common form of empirical processes indexed
by points $x\in\Rd$ and denoted $v_{nF}(x)$, but most of the time we
will be using the function-parametric version of empirical process,
\[
v_{nF}(\phi) = \int_{\Rd} \phi(x)
v_{nF}(\mathrm{d}x) = \frac{1}{\sqrt{n}} \sum_{i=1}^n
\bigl[\phi(X_i) - E\phi(X_i)\bigr], \qquad\phi\in
L_2(F) .
\]
As we know (see, e.g., \cite{vdVW}, Chapter~2), on properly restricted
class of functions $\phi\in{ \Phi}$, the empirical processes
$v_{nF}$ converge to function-parametric Brownian bridge $v_F$. If the
distribution $F$ is uniform on $[0,1]^d$, then $v_F$ becomes a standard
Brownian bridge, which we denote $u$. We recall exact definitions in
the next section.

In this paper, we show that from a certain type of transformation of
$v_{nF}$ a unified approach to distribution free testing of hypothesis
about $F$ is emerging. The approach can be used regardless of whether
the hypothesis is simple or parametric, or whether $F$ is
one-dimensional or multidimensional, and also whether $F$ is continuous
or discrete. The last point is demonstrated in \cite{Khm13} and also
in Corollary~\ref{discrete}, Section~\ref{sec3.1}. We also believe that the
approach is simple to implement: on-going research shows that
parametric families with, multidimensional parameters, as say, family
with 9 parameters, in one of the examples in \cite{Thuong}, can be
studied without noticeable numerical difficulties.

The structure of the transformation in question is the following: let
$K$ be a unitary operator of a certain type, acting on $L_2(F)$, and
consider a transformed process
%
%e1 #&#
%
\begin{equation}
\label{idea} \bigl(K^* v_{nF}\bigr) (\phi) = v_{nF} (K \phi)
.
\end{equation}
The explicit description of the operators we propose to use we defer to
Section~\ref{sec3}, where we show that the processes so obtained will have very
desirable asymptotic properties while being one-to-one transformations
of $v_{nF}$ and, therefore, containing the same amount of ``statistical
information''. As a preliminary illustration of one type of results of
this form, let us formulate the following proposition. It is a
particular case of Theorem~\ref{thm:main1} of Section~\ref{sec3.1}.
%%=======================proposition=============
%
%pr1 #&#

\begin{prop}\label{notag}
Suppose $F$ is an absolutely continuous distribution on $[0,1]^d$
(different from uniform distribution), which has a.e. positive density
$f$. The process $u=\{u(x), x\in[0,1]^d\}$ with the differential
%
%e2 #&#
%
\begin{equation}
\label{unif} u(\mathrm{d}x) = \frac{1}{\sqrt{f(x)}} v_F(\mathrm
{d}x) - \frac{1- \sqrt{f(x)}}{1-
\int_{[0,1]^d} \sqrt{f(y)} \,\mathrm{d}y}
\int_{[0,1]^d}\frac{1}{\sqrt{f(y)}} v_F(\mathrm{d}y) \,\mathrm{d}x
\end{equation}
is the standard Brownian bridge.
\end{prop}

For goodness of fit theory on $\Rd$, this means that with help of a
single stochastic integral above, the asymptotic situation of testing a
simple null hypothesis $F$ can be transformed into the situation of
testing the uniform distribution. In other words, transformation \eqref
{unif} from empirical process $v_{nF}$ possesses the same convenience
for asymptotic statistical inference as the uniform empirical process
in $[0,1]^d$.

As the first step toward \eqref{idea}, in Section~\ref{sec2} below we will find
that there are many more different Brownian bridges than is commonly
realized. We will also see, within the same framework, that although
their distributions remain mutually singular, the boundary between
Brownian bridges and Brownian motion is somewhat blurred and unitary
operators can easily be used to transform Brownian bridges into a
version of bridges, which are ``almost'' Brownian motion. This is
described in Section~\ref{sec3.2}.

%To find unitary operators, such that the transformed processes will be
%asymptotically distribution free was the main motivation for what
%follows. Then one can apply this transformation to empirical processes

%In section 3.1 we describe $U=U_{FG}$ such that if $v_{nF}$ and
%$v_{nG}$ are two empirical processes, based on i.i.d. samples in $\Rd$
%from $F$ and $G$, respectively, then, under general condition, the
%asymptotic distribution of $U v_{nF}$ under $F$ and of $v_{nG}$ under
%$G$ is the same, both converge to $G$-Brownian bridge (see definition
%in section 2). Farther, in Corollary~\ref{cor:g-bridge} we will see,
%that the condition that $F$ is supported on $[0,1]^d$ is, actually,
%unnecessary.

%More importantly, for the approach to work we do not need condition
%that $F$ is continuous, and the same approach can be used for discrete
%distributions.

Let us now briefly outline the situation with distribution free
goodness of fit testing in $\Rd$.

If $F$ is a continuous distribution in $\mathbb R$, and $u_n$ is the
uniform empirical process, then, since \cite{Kol}, we know that
$v_{nF}$ can be transformed to $u_n$ as
\[
v_{nF}(x)=u_n(t), \qquad t=F(x),
\]
or, in function-parametric setting, for $\phi\in L_2([0,1])$
\[
v_{nF}(U\phi)=u_n(\phi),
\]
where $U^*\phi(x) = \phi(F(x))$. It is good to note that this
operator, from $L_2([0,1])$ to $L_2(F)$, is also a unitary operator,
that is,
\[
\int_{y\in\mathbb R} \bigl(U^*\phi\bigr)^2(y) \,\mathrm{d}F(y) =
\int
_{y\in\mathbb R} \phi^2\bigl(F(y)\bigr) \,\mathrm{d}F(y) =\int
_{t\in[0,1]} \phi^2(t) \,\mathrm{d}t ,
\]
although there is little tradition of using this terminology, because
in this situation it looks inconsequential.

An analog of time transformation $t=F(x)$ exists in $\Rd$ as well and
is called the Rosenblatt transformation, \cite{Ros}. In, say,
three-dimensional space, in obvious notation, it has the form
$t_1=F(x_1), t_2=F(x_2|x_1), t_3=F(x_3|x_1,x_2)$. For some reason, and
maybe because dealing with conditional distributions is often awkward,
the transformation is rarely used. It also fails to lead to
distribution free testing, when $F$ depends on a finite-dimensional
parameter (cf. Section~\ref{sec3.3}).

A unitary operator, very different in its nature from time
transformation, was introduced for the empirical processes in
$d$-dimensional time in \cite{Khm88} and \cite{Khm93} and, in two-sample problem, in
\cite{EinKhm}.
%and \cite{Khm93}.
In its origin it is connected with the innovation problem for curves in
Hilbert spaces, \cite{Cra64}, and the theory of innovation
martingales; see, for example, \cite{LipShi}, Section~7.4. In its
simplest form, it is an operator from $\mathcal{L}_F=\{\phi\in
L_2(F)\dvtx \int_{\Rd}\phi(x) \,\mathrm{d}F(x)=0\}$ onto $L_2(F)$
and the result
was a one-to-one transformation from Brownian bridge $v_F$ to Brownian
motion $w_F$. It extends to the case of parametric hypothesis in $\Rd
$. We comment further on it in Section~\ref{sec3.2}.

The approach of this paper seems to us closest to the geometric
argument behind K. Pearson's chi-square statistic, \cite{Pear}; see
also retrospective historic account in \cite{Stig}. The idea itself is
very simple and it is somewhat strange that it was not discovered
before. In the case of one given $F$, the operators involved will map
$\mathcal{L}_F$ in $\mathcal{L}_G$ and subsequently transform one
Brownian bridge, $v_F$, into another Brownian bridge, $v_G$, with $G$
of our choice. Just as Fisher \cite{Fish} and \cite{Fish2}
has extended chi-square theory to the parametric case, our approach as
we said, also extends to the case of parametric families; see Section~\ref{sec3.3}.

Next, in Section~\ref{sec2}, as we said, we present a somewhat broader
definition of Brownian bridges as projected Brownian motions. In
Section~\ref{sec3}, we present the main results. The case of a simple
hypothesis, which also serves as an illustration of the whole approach,
is treated in Section~\ref{sec3.1}, the transformation to ``almost'' Brownian
motion is shown in Section~\ref{sec3.2}, while the case of parametric hypotheses
is considered in Section~\ref{sec3.3}. In Section~\ref{sec3.4}, we discuss the problem of
uniqueness of the proposed transformations. In the last Section~\ref{sec4}, we
illustrate the rate of convergence of transformed empirical processes
to their distribution free limits through the rate of convergence of
the two classical goodness of fit statistics based on these processes:
Kolmogorov--Smirnov statistic and omega-square statistic.

%show the key idea behind transformation \eqref{z} along with its
%different form and an application. In \S3 we describe the case of
%parametric families of discrete distributions. In \S4 we consider the
%extension to the case of continuous time and, in particular, show a
%one-to-one transformation of an $F$-Brownian bridge into a
%$G$-Brownian bridge as well as relatively simple one-to-one
%transformation of $F$-Brownian bridge into a ``slightly restricted''
%standard Brownian motion. In the same section we consider parametric
%families of continuous distributions .

%=================================Sec.2========================

%s2 #&#
\section{Preliminaries: $q$-projected Brownian motions}\label{sec2}

Consider $w_F(\phi), \phi\in L_2(F)$, a function parametric
$F$-Brownian motion, which is a linear functional in $\phi$ and for
each $\phi$ is a Gaussian random variable with mean $0$ and variance
\[
E w_F^2(\phi) =\int_{y\in\Rd}
\phi^2(y) \,\mathrm{d}F(y) = \|\phi\|_F^2 .
\]
This implies that the covariance between $w_F(\phi)$ and $w_F(\tilde
\phi)$ is
\[
E w_F(\phi) w_F(\tilde\phi) =\int_{y\in\Rd}
\phi(y) \tilde\phi(y)\,\mathrm{d}F(y) = \langle\phi, \tilde\phi
\rangle_F .
\]
As far as we are not considering trajectories of $w_F(\phi)$ in $\phi
$, we need only to know that $\phi$ is indeed square integrable with
respect to $F$. For the theory of $w_F(\phi)$ as linear functionals on
$L_2(F)$ and reproducing kernel Hilbert spaces, where they live, we
refer, for example, to \cite{Gross} or the monograph \cite{Kuo}.

Let $v_F(\phi)$ denote the function-parametric $F$-Brownian bridge,
defined as a linear transformation of $w_F$:
%
%e3 #&#
%
\begin{equation}
\label{bvm} v_F(\phi) = w_F(\phi) - \langle\phi,
q_0\rangle_F w_F(q_0).
\end{equation}
Here, we used $q_0$ for the function identically equal to $1$. This
transformation has a particular structure, which is important for what
follows. Namely, we have the following lemma.
%==================== new lemma 1=================
%
%le1 #&#

\begin{lemma}\label{lem:project}
Equality \eqref{bvm} represents $v_F$ as an orthogonal projection of
$w_F$ parallel to the function $q_0$.
\end{lemma}

This statement was initially proved as early as \cite{Khm79}. We show
its proof here for readers' convenience.
\begin{pf*}{Proof of Lemma~\ref{lem:project}}
To shorten notation, denote the right-hand side of
\eqref{bvm} by $\Pi w_F(\phi)$, so that \eqref{bvm} takes the form
$v_F(\phi)=\Pi w_F(\phi), \phi\in L_2(F)$. Then it is easy to see
that
\[
\Pi\Pi w_F = \Pi v_F = v_F,
\]
or $\Pi^2=\Pi$, so that $\Pi$ is indeed a projector. Besides, $\Pi
w_F(q_0)=v_F(q_0)=0$, which, in usual terminology (see, e.g., \cite
{GlaLju}, Section~1.10, and \cite{Kuo}), means that the linear
functional $v_F(\cdot)$ and the function $q_0$ are orthogonal.
\end{pf*}

Substituting the indicator function $\phi= I_{(-\infty, x]}$, from
\eqref{bvm} we obtain
%
%e4 #&#
%
\begin{equation}
\label{Fbridge} v_F(x)= w_F(x) - F(x)w_F(
\infty) ,
\end{equation}
which represents trajectories of $v_F(x)$ as projection of trajectories
of $w_F$. It also leads to the definition of $v_F$ as the Gaussian
process in $x$ with mean $0$ and variance $F(x) - F^2(x)$ (or
covariance $F(\min(x,x')) - F(x) F(x')$).

We can now replace $q_0$ with any other function $q$ of unit
$L_2(F)$-norm. This will lead to the process
%
%e5 #&#
%
\begin{equation}
\label{bvm2} v_F^q(\phi) = w_F(\phi) -
\langle\phi, q\rangle_F w_F(q) ,
\end{equation}
which certainly is again a projection of $w_F$ parallel to $q$ and,
therefore, also could be called Brownian bridge. However, it does not
satisfy the second definition of a bridge. This is more visible in
point-parametric version
%
%e6 #&#
%
\begin{equation}
\label{bvm3} v_F^q(x) = w_F(x) - \int
_{y\leq x} q(y) \,\mathrm{d}F(y) \int_{y\in\Rd}
q(y)w_F(\mathrm{d}y)
\end{equation}
and the variance of $v_F^q$ is of a different form:
%
%e7 #&#
%
\begin{equation}
\label{variance} E\bigl[v_F^q(x)\bigr]^2 = F(x)
- \biggl[\int_{y\leq x} q(y) \,\mathrm{d}F(y)\biggr]^2 ,
\end{equation}
so that if $q\ne q_0$, the second term is not square of the first.
Therefore, even in one-dimensional case, with $F$ being just uniform
distribution on interval $[0,1]$, the distribution of $\max_x
|v_F^q(x)|$ is not Kolmogorov distribution and the distribution of
$\int_0^1[v_F^q(x)]^2 \,\mathrm{d}F(x)$ is not omega-square
distribution unless
$q=q_0$ $F$-a.e. We call $v_F^q(x)$ a slightly longish name of a
\textit{$q$-projected $F$-Brownian motion}. The processes $v_F^q$
arise naturally as weak limits in certain statistical problems and they
will be useful in this paper.

We stress again, that the definition of $v_F^q$ involves two objects --
a distribution $F$ and a function $q\in L_2(F)$.
When $F$ is uniform distribution on $[0,1]^d$ we call $v_F^q$ a
$q$-projected standard Brownian motion (or simply $q$-projected
Brownian motion) and use, most of the time, notation $v^q$ without
index $F$. In the case of general $F$, we would still call $v_F^{q_0}$
a Brownian bridge and often omit $q_0$ from notation. Obviously
$v^{q_0}$ is just a standard Brownian bridge $u$. We formulate the
lemma below for convenience of reference later on.
%========== Lemma~2 (was, initially, lemma 1) ======= new bridge
%===========
%
%le2 #&#

\begin{lemma}\label{lem:newbridge}
Suppose distribution $F$ is supported on the unit cube $[0,1]^d$ and
has a.e. positive density $f$. Suppose $w$ is standard Brownian motion
on $[0,1]^d$ and $v_F$ is defined as in \eqref{bvm} and \eqref
{Fbridge}. Then
\[
\xi(x) = \int_{y\leq x}\frac{1}{\sqrt{f(y)}} v_F(\mathrm{d}y)
\]
is $q$-projected standard Brownian motion with $q=\sqrt{f}$,
%
%e8 #&#
%
\begin{equation}
\label{extra} \xi(x) = w(x) - \int_{y\leq x}\sqrt{f(y)} \,\mathrm
{d}y \int
_{y\in[0,1]^d
}\sqrt{f(y) } w(\mathrm{d}y),
\end{equation}
or, for $\psi\in L_2([0,1]^d)$,
%
%e9 #&#
%
\begin{equation}
\label{extra1} \xi(\psi) = w(\psi) - \langle\psi,\sqrt{f}\rangle
w(\sqrt{f}).
\end{equation}
Conversely, if $\xi$ is $q$-projected standard Brownian motion, then
\[
v_F(x)=\int_{y\leq x} q(y) \xi(\mathrm{d}y)
\]
is $F$-Brownian bridge, as defined in \eqref{Fbridge}, with $F(x)=\int
_{y\leq x} q^2(y) \,\mathrm{d}y$.
\end{lemma}

\begin{pf} The first statement of the lemma follows from the
connection \eqref{Fbridge} between $v_F$ and $w_F$. Indeed, substitute
the normalized differential of $v_F$,
\[
\frac{1}{\sqrt f (y)} v_F(\mathrm{d}y) = \frac{1}{\sqrt f (y)}
w_F(\mathrm{d}y)
- \sqrt f (y) w_F(\Rd) \,\mathrm{d}y,
\]
in the definition of $\xi(x)$ to obtain
\[
\xi(x)=\int_{y\leq x}\frac{1}{\sqrt{f(y)}} w_F(\mathrm{d}y) -\int
_{y\leq
x}\sqrt{f(y)} \,\mathrm{d}y w_F(\infty)
\]
and note that
\[
w(x) = \int_{y\leq x} \frac{1}{\sqrt f (y)} w_F(\mathrm{d}y)
\]
is the standard Brownian motion -- it obviously is $0$-mean Gaussian
process with independent increments and
\[
E\biggl[\int_{y\leq x} \frac{1}{\sqrt f (y)} w_F(\mathrm{d}y)
\biggr]^2= \int_{y\leq x} \frac{1}{f (y)} F(\mathrm{d}y) = x.
\]
Note also that we can write $w_F(\infty)$ as
\[
\int_{y\in[0,1]^d} \sqrt{f(y)} w(\mathrm{d}y).
\]
\upqed
\end{pf}
\begin{Remark*} Note, that the normalization $v_F(\mathrm{d}y)$ by
$\sqrt
f(y)$ does not help to standardize $\xi(\mathrm{d}x)$ -- in \eqref
{extra1} we
still have linear functional $\langle\psi,\sqrt{f}\rangle$, and
thus, the dependence on $F$ in $\xi$ is still present. This was well
understood for a very long time, and it is quite unexpected that using
one extra stochastic integral (see Proposition~\ref{notag}), the
standardization becomes possible.
\end{Remark*}

The normalization by $1/\sqrt{f}$ used in the lemma is a particular
form of the more general mapping.
Namely, let $G$ be another distribution on $\Rd$, which is absolutely
continuous with respect to $F$. Then the function
%
%e10 #&#
%
\begin{equation}
\label{l} l(x)=\sqrt{\frac{\mathrm{d}G}{\mathrm{d}F} (x)}
\end{equation}
belongs to $L_2(F)$. and moreover, if $\psi\in L_2(G)$, then $l\psi
\in L_2(F)$ and $\| \psi\|_G=\|l\psi\|_F$. If distributions $G$ and
$F$ are equivalent (mutually absolutely continuous), then the inverse
is also true: if $\phi\in L_2(F)$ then $\phi/l \in L_2(G)$, and the
norm is preserved. This, in particular, means that re-normalization of
$F$-Brownian motion into $G$-Brownian motion is straightforward: if
$w_F(\phi)$ is an $F$-Brownian motion in $\phi\in L_2(F)$, then
$w_F(l\psi) = w_G(\psi)$ is a $G$-Brownian motion in $\psi\in
L_2(G)$. This, we repeat, does not extend to $v_F^q$ and $v_F$ -- the
distribution of, say, $v_F(l\psi)$ depends on both $F$ and $G$. The
first theorem in Section~\ref{sec3} below shows, however, that a simple
isomorphism exists.

To describe one more object we consider in this paper, complement $q_0$
by a sequence of orthonormal functions $q_1,\dots, q_\kappa$, which
are also orthogonal to $q_0$, and consider the process\vspace*{1.5pt}
\[
\hat{ v}_F(\phi) = w_F(\phi) - \sum
_{i=0}^\kappa\langle q_i,\phi
\rangle_F w_F(q_i) .
\]
Similar to what we said about $v_F$, the process $\hat{ v}_F$ is the
orthogonal projection of $w_F$ parallel to the functions $q_0,\dots,
q_\kappa$. We still call $\hat{ v}_F$ a $q$-projected $F$-Brownian
motion. It may be that notation $\hat{v}^q_F$ is used again, but when
$q$ is a vector function, there is no other ``more traditional'' notion
to be confused with $\hat{v}_F$; so we skip $q$ as an upper index.

The role of the process $\hat{v}_F$ becomes clear when we examine
asymptotic behavior of the parametric empirical process. Consider the
problem of testing parametric hypothesis that the distribution function
of $X_i$s belongs to a given family of distribution functions $F_\theta
(x)$, depending on a finite-dimensional parameter $\theta$. The value
of this parameter is not prescribed by the hypothesis and has to be
estimated using the sample $X_1,\dots, X_n$. Denote\vspace*{1.5pt}
\[
v_{n}(B, \hat{\theta}_n)=\sqrt{n} \bigl[F_n(B)
- F_{\hat{\theta}_n} (B)\bigr]
\]
the parametric empirical process (indexed by sets). (Note that, in
presence of $\theta$ and $\hat{\theta}$, one can skip index $F$ in
notation.) As has been known since Kac \textit{et al}. \cite{KKW} and
later Durbin \cite{Dur73}
and other work, the asymptotic behavior of empirical processes with
estimated parameters is different from that of $v_{nF}$, and in
particular, its limit distribution depends not only on the true value
of the parameter but also on the score function. However, we can say
more.

Namely, under usual and mild assumptions (see, e.g., \cite{Cra46},
Chapter~3, and see the modern exposition in \cite{vdV}, Section~5),
the MLE $\hat{\theta}_n$ possesses an asymptotic representation\vspace*{1.5pt}
\[
\sqrt n (\hat{\theta}_n -\theta)=\Gamma^{-1}_{\theta}
\int_{x\in
\Rd} \frac{\dot f_\theta(x)}{f_\theta(x)} v_{nF}(\mathrm{d}x,
\theta) +
\mathrm{o}_P(1) ,\qquad n\to\infty,
\]
where we denote by $f_\theta$ and $\dot f_\theta$ the hypothetical
density and the vector of its derivative with respect to parameter
$\theta$ and denote\vspace*{1.5pt}
\[
\Gamma_\theta= \int_{y\in\Rd}\frac{\dot f_\theta(y) \dot
f_\theta^T(y)} {f_\theta(y)} \,\mathrm{d}y
\]
the Fisher information matrix. Consequently, the parametric empirical
process has asymptotic expansion
%
%e11 #&#
%
\begin{eqnarray}
\label{vhat} v_n(B, \hat{\theta}_n) &=&
v_n(B, \theta) - \int_{B} \frac{\dot
f_\theta^T (y)}{f_\theta(y)}
F_{\theta} (\mathrm{d}y) \Gamma^{-1}_\theta\int
_{y\in\Rd} \frac{\dot f_\theta(y)}{f_\theta(y)} v_n(\mathrm{d}y,
\theta)
+\mathrm{o}_P(1) ,
\nonumber
\\[-8pt]\\[-8pt]
&=& v_n(B, \theta) - \int_{B}
\beta_F^T(y) F_{\theta} (\mathrm{d}y) \int
_{y\in\Rd} \beta_F(y)v_n(\mathrm{d}y,\theta) +
\mathrm{o}_P(1),\nonumber
\end{eqnarray}
where
\[
\beta_F(x) = \Gamma^{-1/2}_\theta\frac{\dot f_\theta(x)}{f_\theta
(x)}
.
\]
As shown in \cite{Khm79} (see also \cite{KhmKo04}, Section~2.2),
the main part of this expansion represents $v_n( \cdot, \hat{\theta
}_n)$ as the orthogonal projection of $v_n( \cdot,\theta)$ parallel
to the normalized score function $\beta_F$ and, therefore, the limit
in distribution of $v_n(\cdot, \hat{\theta}_n)$ can be written (in
function-parametric form) as
\[
v_F(\phi) - \bigl\langle\beta_F^T,\phi
\bigr\rangle_F v_F(\beta_F).
\]
At the same time, the score function $\beta_F$ is orthogonal to the
function $q_0$ and its coordinates are orthonormal and will play the
role of functions $q_1,\dots, q_\kappa$ above. Therefore,
substituting representation \eqref{bvm} of $v_F$ through $w_F$, we see
that the limit in distribution of the process $v_N(\cdot,\hat{\theta
}_n)$ is $\hat{v}_F$.

It is well known that the actual weak convergence statement in
function-parametric set-up requires some restriction on the underlying
class of functions $\phi$, but these restrictions are well understood
and we refer readers to \cite{vdVW}. For an earlier proof in Skorohod
space, see \cite{Koul69} and \cite{Dur73}, while for the proof in
$L_2(F)$ see \cite{Khm79}.

%==========================Sec 3==============================

%s3 #&#
\section{The main result and its corollaries}\label{sec3}

The main geometric idea in this paper can be described as follows. When
testing for fixed distribution $F$, the corresponding empirical
processes will converge to $v_F$, which is an orthogonal projection of
the Brownian motion. When testing for a different $G$, there will be
convergence to $v_G$, which is also an orthogonal projections of
Brownian motion. However, we will see that if $G$ and $F$ are
equivalent; these projections can be ``rotated'' to each other. The
unitary operators involved in this rotation form a group, transient on
the class of all Brownian bridges with all $G$ equivalent to the $F$.
In other words, the problem of testing $F$ can be mapped to the problem
of testing $G$ and vice versa, and these, seemingly distinct problems
are not distinct problems, but form one equivalence class. Therefore,
one representative of each equivalence class is sufficient, and we
propose a form of such representative. Since the processes $\hat{v}_F$
and $\hat{v}_G$ are both orthogonal projections as well, the idea of
unitary transformation extends to the parametric classes of
distributions.

%======================Simple hypothesis=============

%s3.1 #&#
\subsection{The case of fixed $F$}\label{sec3.1}
Although the following Theorem~\ref{thm:main1} is generalized by
Theorem~\ref{thm:main2}, by starting with the case of one fixed $F$
and giving an independent proof we hope to make the overall
presentation more transparent.

Consider an operator on $L_2(F)$\vspace*{1pt}
%
%e12 #&#
%
\begin{equation}
\label{unitary} K=I -\frac{2}{\|l-q_0\|_F^2} (l-q_0)\langle
l-q_0,\cdot\rangle_F,
\end{equation}
where $I$ is identity operator and $l$ is the function defined in
\eqref{l}, while the function $q_0$ identically equals 1. Below we
will also need the linear subspace $\mathcal{L}=\mathcal{L}(q_0,l)$,
generated by functions $q_0$ and $l$ and functions $l_\perp$ and
$q_{0,\perp}$, which are parts of $l$ and $q_0$, orthogonal to $q_0$
and $l$, respectively,\vspace*{1pt}
\[
l_\perp=l-\langle l,q_0\rangle_F
q_0,\qquad q_{0\perp}=q_0-\langle l,q_0
\rangle_F l .
\]
It is clear that\vspace*{1pt}
\[
\biggl(q_0,\frac{1}{\|l_\perp\|} l\biggr) \quad\mbox{and}\quad
\biggl(l,
\frac{1}{\|q_{0\perp
}\|} q_{0\perp}\biggr)
\]
form two orthonormal bases of $\mathcal{ L}$.

The operator $K$ has the following properties.
%=====================lemma 3===============
%
%le3 #&#

\begin{lemma}\label{lem:Koperat} \textup{(i)} Operator $K$ is a (self-adjoint)
unitary operator on $L_2(F)$, $\|K\phi\|_F=\|\phi\|_F$, such that\vspace*{1pt}
\[
K\phi= \phi,\qquad \mathit{if} \phi\perp l, q_0 \quad\mbox
{and}\quad Kl=q_0,\qquad
Kq_{0\perp} = l_\perp,\qquad  \mathit{while}\ Kq_0=l.
\]

\textup{(ii)} Coordinate representation of this operator is\vspace*{1pt}
\[
K=I_{\mathcal{L}\perp} + q_0\langle l,\cdot\rangle_F +
l_{\perp} \langle q_{0\perp},\cdot\rangle_F,
\]
where $I_{\mathcal{L}\perp}$ is the projection operator on the
subspace of $L_2(F)$ orthogonal to $\mathcal{ L}$.
\end{lemma}

The reader can easily verify the lemma. Part (i) is needed just
below, part (ii) will be useful to draw similarity with Section~\ref{sec3.3}.
below. Note that one could use a similar unitary operator, with $l-q_0$
replaced by $l+q_0$. We chose the present form only because the norm $\|
l-q_0\|_F$ is a well-known object -- the Hellinger distance between
distributions $F$ and $G$. To what extent the choice of $K$ is unique
is discussed in Section~\ref{sec3.4}. Note also that\vspace*{1pt}
%
%e13 #&#
%
\begin{equation}
\label{hellinger1} \|l-q_0\|_F^2 = 2\int
_{y\in\Rd} \bigl(1- l(y)\bigr)\,\mathrm{d}F(y)=-2\langle l-q_0,
q_0\rangle_F .
\end{equation}
%
%=======================Theorem~2 between bridges=========

%th2 #&#
%
\begin{them}\label{thm:main1}
Suppose distribution $G$ is absolutely continuous with respect to
distribution $F$ (and different from $F$). If $v_F$ is $F$-Brownian
bridge, then the process with differential
%
%e14 #&#
%
\begin{eqnarray}
\label{G-to-F} v_G(\mathrm{d}x) &= & l(x) v_F(\mathrm{d}x) \nonumber
\\[-8pt]\\[-8pt]
&&{}- \int
_{y\in\Rd} l(y) v_F(\mathrm{d}y) \frac{1}{ 1- \int_{y\in\Rd}
l(y) \,\mathrm{d}F(y)} \bigl[l^2(x) - l(x)\bigr] f(x) \,\mathrm{d}x
\nonumber
\end{eqnarray}
is $G$-Brownian bridge.

If distributions $G$ and $F$ are equivalent, that is, if $l=\sqrt
{\mathrm{d}G/\mathrm{d}F}$ is positive $F$-a.e., then \eqref{G-to-F}
is one-to-one.
\end{them}

If $F$ is an absolutely continuous distribution on the unit cube
$[0,1]^d$ and its density $f$ is positive a.e., while $G$ is uniform on
this cube, then $l(x)=1/\sqrt{f(x)}$ and we obtain the transformation
of $F$-Brownian bridge into the standard Brownian bridge, already given
in Proposition~\ref{notag}.
\begin{Remark*} It was interesting to realize that $v_G$ in \eqref
{G-to-F} remains $G$-Brownian bridge even if $\mathrm{d}G/\mathrm
{d}F$ can be $0$ on a
set of positive probability $F$.
\end{Remark*}
\begin{pf*}{Proof of Theorem~\ref{thm:main1}} As we know, for any
function $\psi\in L_2(G)$, under
our conditions, $l\psi\in L_2(F)$. Since $v_F(l-q_0)=v_F(l)$, the
function-parametric form of \eqref{G-to-F} is
\[
v_G(\psi)=v_F(K \phi) , \qquad\mbox{with } \phi= l\psi.
\]
We need to show that the covariance operator of $v_G$ is that of
$G$-Brownian bridge. For this it is sufficient to consider the variance
of $v_G (\psi)$,
\[
E\bigl[v_G(\psi)\bigr]^2 = E\bigl[v_F(K
\phi)\bigr]^2 = \|K \phi\|^2_F - \bigl[
\langle K \phi, q_0\rangle_F\bigr]^2 .
\]
However,
\[
\|K \phi\|^2_F=\|\phi\|^2_F =
\|\psi\|^2_G ,
\]
and, using \eqref{hellinger1}, we obtain
\begin{eqnarray*}
\langle K \phi, q_0\rangle_F &= &\langle\phi,
q_o\rangle_F - \frac
{2}{\|l-q_0\|_F^2} \langle
l-q_0, q_0\rangle_F\langle
l-q_0,\phi\rangle_F =\langle\phi, l
\rangle_F
\\
& =& \langle l \psi, l \rangle_F = \langle\psi, q_0
\rangle_G .
\end{eqnarray*}
Therefore,
\[
E\bigl[v_G(\psi)\bigr]^2 = \| \psi\|^2_G
- \bigl[\langle\psi, q_0\rangle_G\bigr]^2 ,
\]
which is the expression for the variance of $G$-Brownian motion.
\end{pf*}

Although any distribution $F$ in $\Rd$ can be mapped to a distribution
on the unit cube, in some cases this mapping may involve unpleasant
technicalities. Corollary~\ref{cor:g-bridge} helps to make this
mapping very simple, and actually unnecessary, in a wide class of situations.
The idea is that $v_F$ can be transformed into $v_G$, and for this $G$
the mapping to the unit cube will be immediate. Namely, choose $d$
densities $g_1,\dots,g_d$ on $\mathbb R$, and let
\[
\label{margin} g(x)=\prod_{i=1}^d g_i(x_i)
.
\]
Denote $t_i=\int_{-\infty}^{x_i} g_i(s) \,\mathrm{d}s, i=1,\dots,
d$. Then
%
%e15 #&#
%
\begin{equation}
\label{time} \prod_{i=1}^d t_i =
\prod_{i=1}^d \int_{-\infty}^{x_i}
g_i(s) \,\mathrm{d}s
\end{equation}
is direct $d$-dimensional analogue of Kolmogorov time transformation
$t=G(x)$ on the real line.
It seems clearer to give the formulation of the next statement for
rectangles rather than for general Borel sets $B$.
%
%========================= Cor ==Time transformation ===============
%
%co3 #&#

\begin{cor}\label{cor:g-bridge} Suppose $g_1,\dots,g_d$ are such that
the distribution $G$ with density
$g$ is absolutely continuous with respect to $F$. Suppose the points
$t\in[0,1]^d$ and $x\in\Rd$ are connected as in \eqref{time}. If
$v_F$ is $F$-Brownian bridge and $v_G$ is its transformation \eqref
{G-to-F}, then the process $u$,
\[
\label{unit-d} u(t) = v_G(x),
\]
is a standard Brownian bridge on $[0,1]^d$.
\end{cor}

It is now clear that there is no need to perform the time
transformation \eqref{time}, because it is obvious how to choose test
statistics from $v_G$, which are invariant under this transformation.
For example, for $G$ as in \eqref{time}, the statistics
%
%e16 #&#
%
\begin{equation}
\label{df-stat} \sup_{x\in\Rd} \bigl|v_G(x)\bigr| \quad\mbox{and}\quad
\int_{x \in\Rd} v_G^2(x) \,\mathrm{d}G(x)
\end{equation}
have distributions independent from $G$ and, hence, from the initial
distribution $F$.
On the other hand, the class of distributions $F$ for which the product
distribution $G$ exists, and then there are infinitely many of them, is
broad: any distribution which has rectangular support, whether bounded
or unbounded, is such a distribution. Equivalently, if the copula
function corresponding to $F$ has positive density on $[0,1]^d$, then
$G$ exists (see, e.g., \cite{Nel} and \cite{Joe} for such examples)
and one can choose $g_1,\dots, g_d$ as marginal densities of $F$.

As an immediate consequence of Theorem~\ref{thm:main1} for finite $n$,
we have the following weak convergence statement. Consider the process
%
%e17 #&#
%
\begin{eqnarray}
\label{G2F-n} \tilde v_n(x) &= &\int_{y\leq x} l(y)
v_{nF}(\mathrm{d}y)\nonumber\\[-8pt]\\[-8pt]
&&{}- \int_{y\in\Rd} l(y) v_{nF}(\mathrm{d}y)
\frac{1}{ 1- \int_{y\in\Rd} l(y) \,\mathrm{d}F(y)} \biggl[G(x) -
\int_{y\leq x} l(y) \,\mathrm{d}F(y)
\biggr].
\nonumber
\end{eqnarray}
%
%=====================corollary==finite n=============
%
%co4 #&#

\begin{cor}\label{cor:finite n} Let $v_G$ be the point-parametric
$G$-Brownian bridge defined in \eqref{G-to-F}. Then, as $n\to\infty$,
\[
\tilde v_{n} \stackrel{\mathcal{D}(F)} {\longrightarrow}
v_G .
\]
\end{cor}

In other words, the limit distribution of $\tilde v_{n}$ under $F$ is
the same as the limit distribution of empirical process $v_{nG}$ under
$G$. If $F$ has a rectangular support, then, as noted above, $G$ of
product form exists. Then, using \eqref{time}, $\tilde v_{n}$ can
further be transformed into a process, which under $F$, converges in
distribution to the standard Brownian bridge $u$. In other words,
construction of asymptotically distribution free test statistics from
$\tilde v_{n}$ becomes obvious, cf. \eqref{df-stat}.

For the proof of this corollary, note that the weak convergence
statement for the first integral in \eqref{G2F-n} as a process in $x$
easily follows from, say Theorem~2.5.2 of \cite{vdVW}, as it can be
viewed as statement for function-parametric process indexed by
functions $l(y)\mathbf{1}_{(-\infty, x]}(y)$, which certainly satisfy
the conditions of that theorem. Convergence of the second integral,
with respect to $v_{nF}$, is also clear, while the rest is a fixed
deterministic function.

Our last corollary in this section uses the fact that in Theorem~\ref
{thm:main1} we did not need absolute continuity of $G$ and $F$ with
respect to Lebesgue measure, but only absolute continuity of $G$ with
respect to $F$. Therefore, we can consider discrete distributions with
infinitely many positive probabilities.

Suppose $\mathcal{ X}$ is a countable collection of points of, say,
$\Rd$, and $F$ is a (discrete) probability distribution on $\mathcal
{X}$ with probabilities $p(x)>0$. Suppose $G$ is another distribution
on $\mathcal{ X}$ with probabilities $\pi(x)$. Definition of $v_F$
and $w_F$, as Gaussian processes with prescribed covariance, carries
out to the case of discrete $F$ without change. The differential
$v_F(\mathrm{d}x)$ will now be a jump of $v_F$ at $x\in\mathcal{ X}$
and will
be $0$ at any other $x$. Thus, we obtain the following statement. It
can be viewed as an extension of Theorem~1, (ii), of \cite{Khm13} for
$m=\infty$. In its form, it is no different from \eqref{G-to-F} but
for the fact that $F$ discrete.
%
%==============================Cor discrete ===========
%
%co5 #&#

\begin{cor}\label{discrete} For $x\in\mathcal{X}$, let $l(x)=\sqrt
{\pi(x)/p(x)}$. If $v_F$ is $F$-Brownian bridge (on $\mathcal{X}$),
then the process
%
%e18 #&#
%
\begin{eqnarray}
\label{disc} v_G(B)& = &\int_{y\in B\cap\mathcal{X}} l(y)
v_F(\mathrm{d}y)\nonumber\\[-8pt]\\[-8pt]
&&{}- \int_{y\in
\mathcal{X}} l(y) v_F(\mathrm{d}y)
\frac{1}{ 1-\int_{y\in\mathcal{ X}} l(y) \,\mathrm{d}F(y)} \int
_{y\in B\cap
\mathcal{X}} \bigl(\mathrm{d}G(y) - l(y) \,\mathrm{d}F(y)\bigr)
\nonumber
\end{eqnarray}
is $G$-Brownian bridge.
\end{cor}

Weak convergence statement in discrete case is very simple: with no
possibility of misunderstanding, denote the transformation in \eqref
{disc} applied to $v_{nF}$ again by $\tilde v_{n}$. For any functional,
or statistic, $S(\tilde v_n)$ based on this $\tilde v_n$, which has the
property that for arbitrary small $\epsilon>0$ there is a finite
collection of points $\mathcal{X}_\epsilon$, and a functional
$S_\epsilon(\tilde v_n)$, which depends only on $\tilde v_n(x), x\in
\mathcal{X}_\epsilon$ and is such that
\[
{\mathbb P} \bigl(\bigl|S(\tilde v_n) - S_\epsilon(\tilde
v_n)\bigr| > \epsilon\bigr)< \epsilon,
\]
for all sufficiently large $n$, then
\[
S(\tilde v_n) \stackrel{d(F)} {\longrightarrow} S(\tilde
v_G) .
\]

%s3.2 #&#
\subsection{Mapping to Brownian motion}\label{sec3.2}
Consider one more form of unitary transformations, applied to a
Brownian bridge. It takes somewhat unusual form and seems important in
its own right. In particular, it shows how blurred the difference
between Brownian motions and Brownian bridges can become and, using
some freedom of speech, that ``a Brownian motion can be also a Brownian
bridge''.

Let the distribution $F$ be supported on $[0,1]^d$ and have there a.e.
positive density $f$. Let $A$ be a fixed subset of $[0,1]^d$ and let
$\eta_A(x)$ denote the square root of a density concentrated on $A$,
that is, $\eta_A(x) = 0$ if $ x\notin A$, and $\int_A\eta_A^2(x)\,
\mathrm{d}x =1$.
It is appropriate to think about $A$ as a ``small'' set, although there
will be no formal requirement on this. As we know (see Lemma \ref
{lem:newbridge} in Section~\ref{sec2}), the process with the differential
\[
\xi(\mathrm{d}x)=\frac{v_F(\mathrm{d}x)}{\sqrt{f(x)}}
\]
is the $\sqrt f$-projected standard Brownian motion. At the same time,
the process with the differential
%
%e19 #&#
%
\begin{equation}
\label{db} b(\mathrm{d}x) = w(\mathrm{d}x) - \eta_A(x)\int_A
\eta_A(y) w(\mathrm{d}y) \,\mathrm{d}x
\end{equation}
is the $\eta_A$-projected standard Brownian motion, cf. \eqref
{extra}, and satisfies orthogonality condition
\[
\int_A \eta_A(x) b(\mathrm{d}x) =0 .
\]
In other words, the distribution $F$ is, in the both cases, just
uniform distribution on $[0,1]^d$, but the processes are projected
parallel to different functions. What we want to do now is to rotate
$\xi$ to $b$.

Since $\xi(\phi)$ and $b(\phi)$ are now defined on $L_2([0,1]^d)$,
for our rotation we need to use operator
\[
U^* = I - \frac{2}{\|\eta_A - \sqrt{f}\|^2}(\eta_A - \sqrt
{f})\langle
\eta_A-\sqrt{f},\cdot\rangle,
\]
which is (self-adjoint) unitary operator on $L_2([0,1]^d)$ and maps
$\eta_A$ to $\sqrt{f}$. (Here and in the proof below, in inner
products and norms in $L_2([0,1]^d)$ we skip the index $F$.) The
result, the process
\[
\xi\bigl(U^*\phi\bigr) ,
\]
is what we consider in the next statement. Although a general principle
here remains the same, we believe it is more convenient to formulate it
as a theorem and give the proof.
%
%====================Corollary == Brownian motion ===========
%
%th6 #&#

\begin{them}\label{cor:nine}
Choose $\eta^2_A$ to be a density on $A$. With above assumption on
$F$, the process with differential\vspace*{1pt}
\begin{eqnarray*}
b(\mathrm{d}x)& =& \frac{v_F(\mathrm{d}x)}{\sqrt{f(x)}}\\
&&{} - \int_{y\in A}
\eta_A(y) \frac{v_F(\mathrm{d}y)}{\sqrt{f(y)}} \frac{1}{1 - \int
_{y\in A} \eta_A(y)
\sqrt{f(y)}\,\mathrm{d}y} \bigl(\eta_A(x) -\sqrt{f(x)}\bigr) \,
\mathrm{d}x
\end{eqnarray*}
is a standard Brownian motion on $[0,1]^d\setminus A$, while
\[
\int_{y\in A}\eta_A(y) b(\mathrm{d}y) =0 .
\]
\end{them}

In other words, $b$ is $\eta_A$-projected standard Brownian motion.
\begin{pf*}{Proof of Theorem~\ref{cor:nine}}
The last equality follows from definition of $b$.
Using the process $\xi$, see Lemma~\ref{lem:newbridge}, we easily see
that the function-parametric form
of the process $b$ is
\[
b(\phi) = \xi(\phi) - \frac{2}{\|\eta_A - \sqrt{f}\|^2} \xi(\eta
_A) \langle\phi,
\eta_A - \sqrt{f} \rangle.
\]
Now note that from the definition of $\xi$ in \eqref{extra1} it
easily follows that for $\phi, \tilde\phi\in L_2([0,1]^d)$
\[
E\xi(\phi) \xi(\tilde\phi) = \langle\phi,\tilde\phi\rangle-
\langle\phi,\sqrt{f}
\rangle\langle\tilde\phi,\sqrt{f}\rangle.
\]
Also note that
\[
\|\eta_A - \sqrt{f}\|^2= 2 \bigl(1- \langle
\eta_A,\sqrt{f}\rangle\bigr),
\]
which will somewhat simplify notations below. Thus, we obtain
\begin{eqnarray*}
E b(\phi)^2 &= &\langle\phi,\phi\rangle- \langle\phi,\sqrt{f}
\rangle^2
\\
&&{}- \frac{2}{1 - \langle\eta_A, \sqrt{f}\rangle} \bigl( \langle
\phi, \eta_A \rangle- \langle
\phi,\sqrt{f}\rangle\langle\eta_A, \sqrt{f}\rangle\bigr)
\langle
\phi,\eta_A - \sqrt{f} \rangle
\\
&&{}+ \frac{1}{(1 - \langle\eta_A, \sqrt{f}\rangle)^2}\bigl(1-
\langle\eta_A,\sqrt{f}
\rangle^2\bigr) \langle\phi, \eta_A - \sqrt{f}\rangle
^2 ,
\end{eqnarray*}
or
\begin{eqnarray*}
E b(\phi)^2 &= &\langle\phi,\phi\rangle- \langle\phi,\sqrt{f}
\rangle^2 \\
&&{}- \frac{\langle\phi,\eta_A - \sqrt{f} \rangle}{1 -
\langle\eta_A, \sqrt{f}\rangle} \bigl[ 2 \langle\phi,
\eta_A \rangle- 2 \langle\phi,\sqrt{f}\rangle\langle
\eta_A, \sqrt{f}\rangle
- \bigl(1+ \langle\eta_A,\sqrt{f}\rangle\bigr) \langle\phi,
\eta_A - \sqrt{f}\rangle\bigr] ,
\end{eqnarray*}
and after simplifications within the square brackets we finally obtain
\begin{eqnarray*}
E b(\phi)^2& = &\langle\phi,\phi\rangle- \langle\phi,\sqrt{f}
\rangle^2 - \langle\phi,\eta_A - \sqrt{f} \rangle\langle
\phi, \eta_A + \sqrt{f}\rangle
\\
&=& \langle\phi,\phi\rangle- \langle\phi,\eta_A\rangle^2
,
\end{eqnarray*}
which proves the claim: restriction $x\in[0,1]^d\setminus A$ is
equivalent to restriction that $\phi$ is orthogonal to all $\eta_A$
with given $A$, in which case we obtain the variance of just Brownian
motion, while if $\phi=\eta_A$ the variance of $b(\eta)$ is $0$.
\end{pf*}
\begin{Remark*} If we choose $\eta_A^2$ as the uniform density on
$A$, $\eta_A(x)= I_A(x)/\Delta$, with $\Delta=\mu_d(A)$, then the
process $b$, or rather the finite $n$-version of the process, is certainly
%
%e20 #&#
%
\begin{equation}
\label{bm21} b_n(\mathrm{d}x)= \frac{v_{nF}(\mathrm{d}x)}{\sqrt
{f(x)}} - \int
_A \frac
{v_{nF}(\mathrm{d}y)}{\sqrt{f(y)}} \frac{1}{\sqrt{\Delta} - \int
_A \sqrt
{f(y)}\,\mathrm{d}y} \bigl(
\eta_A(x) -\sqrt{f(x)}\bigr) \,\mathrm{d}x
\end{equation}
which integrates to $0$ on $A$. This, however, should not be perceived
as a ``loss of observations on $A$'': the integral with respect to
$v_{nF}$ over $A$ enters the differential of $b_n$ at all $x\notin A$.
\end{Remark*}
\begin{Remark*} If we choose $\eta_A(x)^2$ as the conditional
density of $F$ given $A$,
$\eta_A^2(x)=1_A(x) f(x) /\allowbreak  F(A)$, then
%
%e21 #&#
%
\begin{equation}
\label{bm2} b_n(\mathrm{d}x) = \frac{v_{nF}(\mathrm{d}x)}{\sqrt
{f(x)}} +v_{nF}(A)
\frac{1}{\sqrt{
F(A)} - F(A)} \sqrt{f(x)} \,\mathrm{d}x ,\qquad x\notin A,
\end{equation}
is another asymptotically Brownian motion on $[0,1]^d\setminus A$. In
this version integration over $A$, where $f$ may happen to be
numerically small, is replaced by $v_n(A)$.
The latter is simpler to calculate and may have better convergence
properties than the integral $\int_A (1/\sqrt{f(y)}) v_{nF} (\mathrm{d}y)$
(cf. Figures~\ref{fig:normedmotion} and \ref{fig:motion} of
Section~\ref{sec4}).
\end{Remark*}

A one-to-one transformation, of a different nature, of a Brownian
bridge to a Brownian motion was earlier suggested in \cite{Khm88}
and \cite{Khm93}. It is interesting to compare that transformation
with the present one. For this we need a so called scanning family of
subsets $S_t,0\leq t\leq1$, of $[0,1]^d$, which is increasing,
$S_t\subseteq S_{t'}$ for $t<t'$, and such that $\mu_d(S_0)=0, \mu
_d(S_1)=1$ and $\mu_d(S_t)$ is continuous in $t$. Then, with $\xi$ as
above, the process
\[
{\tilde b}(C, t)= \xi(C\cap S_t) - \int_0^t
\frac{ \int_{S_\tau^c}
\sqrt{f(y)} \xi(\mathrm{d}y)}{1 - F(S_\tau)} d\int_{C\cap S_\tau}
\sqrt{f(z)} \,\mathrm{d}z
\]
is not only a Gaussian martingale in $t$, but also has independent
increments in $C\subseteq S_t$, so that $\tilde b(C,1)$ is a Brownian
motion in $C\subseteq[0,1]^d$. The latter expression is a
multidimensional extension of the classical situation for $d=1$ and
$f=1$ on $[0,1]$, when the $\xi(t)=u(t)$ is the standard Brownian
bridge. Indeed, from the above we obtain the well-known representation
of $u(t)$ as a Gaussian semimartingale
\[
\tilde b(\mathrm{d}t) = u(\mathrm{d}t) + \frac{u(t)}{1-t} \,\mathrm
{d}t ,
\]
where $\tilde b$ is Brownian motion. In statistical context, see its
use in \cite{ABGK} and \cite{Khm81}, Section~1; see also \cite{ShW},
Chapter~6. The inverse of this representation,
\[
u(t)=(1-t)\int_0^t\frac{\tilde b(\mathrm{d}s)}{1-s} ,
\]
was used for statistical purpose as early as \cite{Doo49}.

The transformation of Theorem~\ref{cor:nine} is simpler; for
$A=[0,\Delta]\subseteq[0,1]$ it takes the form
\begin{eqnarray*}
b(\mathrm{d}t) &=& u(\mathrm{d}t) - \frac{u(\Delta)}{\Delta} \,
\mathrm{d}t ,\qquad t\leq\Delta,
\\
b(\mathrm{d}t) &=&u(\mathrm{d}t) - \frac{u(\Delta)}{\sqrt{\Delta}
- \Delta} \,\mathrm{d}t, \qquad t>\Delta,
\end{eqnarray*}
and represents a Brownian bridge on $[0,\Delta]$ and Brownian motion
on $[\Delta,1]$. Although in last three displays the same process $u$
is transformed and the same, in distribution, process is obtained on
$[\Delta,1]$ as a result, the transformations are very different.

%
%Some heuristic comments may be helpful here, although they do not
%explain the following corollary fully. Let $e_1,\dots, e_n,\dots$, be
%an orthonormal sequence in $L_2([0,1]^d)$. Choose $e_m$ with large $m$
%and fix it. Returning for a moment to the function-parametric
%notations, the process $\xi_1(\phi)=w(\phi) - w(e_1)\langle e_1,\phi
%\rangle$ is $e_1$-orthogonal Brownian bridge, while $\xi_m(\phi)=w(
%\phi) - w(e_m)\langle e_m,\phi\rangle$ is $e_m$-orthogonal Brownian
%bride. Orthogonality to $e_1$ is an ``essential'' restriction and
%trajectories of $w$ and $\xi_1$ will ``look different'', but
%orthogonality to $e_m$ can be much ``weaker'' restriction, much less
%felt, because for ``most'' functions the Fourier coefficient $\langle
%e_m,\phi\rangle$ will be small. Yet, there is a simple isometry
%between the two processes:
%$$ \xi_m(\phi)\stackrel{d}{=}\xi_1(\phi) - \xi_1(e_m)(\langle e_m,\phi
%\rangle- \langle e_1,\phi\rangle). $$ As the next statement shows,
%$e_1$ and $e_m$ do not have to be orthogonal, but only linearly
%independent.

%===================section ======== parametric case=========
%s3.3 #&#
\subsection{Parametric family of distributions}\label{sec3.3}
We extend now the results for the case of fixed $F$ to the parametric
case. Namely, along with distribution $F_\theta$ and its orthonormal
score function $\beta_F$, consider now another distribution $G$
together with orthonormal vector $(r_1,\dots, r_\kappa)^T$, with
coordinates in $L_2(G)$, of the same dimension as $\beta_F$. One may
think about this vector $\beta_G$ as a score function of a more or
less fictitious parametric family to which $G$ belongs. Let us augment
both score functions by a function identically equal~1. If $G$ is
absolutely continuous with respect to $F$, then the vector
$(l,\mathit{lr}_1,\dots,\mathit{lr}_d)$ is orthonormal in $L_2(F)$.

Use notation $\hat{\mathcal{L}}$ for a subspace of functions
\[
\hat{\mathcal{L}}= \mathcal{L} (q_0,\dots, q_\kappa, l,
\dots, \mathit{lr}_\kappa)\subset L_2(F),
\]
where we recall, $q_0=1$ and $q_i, i=1,\dots,\kappa$, are coordinate
functions of $\beta_F$. In the subspace $\hat{\mathcal{L}}$,
consider two bases. One, the $a$ basis, has coordinate functions
$a_i=q_i$ for $i\leq\kappa$ while $a_i, i=\kappa+1, \dots, 2\kappa
+1$, is any orthonormal sequence, which complements $a_0,\dots,
a_\kappa$ to a basis in $\hat{\mathcal{ L}}$. The other, $b$ basis,
has coordinate functions $b_i=l r_i, i\leq\kappa$, and $b_i, i=\kappa
+1, \dots, 2\kappa+1$, can be any orthonormal sequence, which
complements $b_0,\dots, b_\kappa$ to a basis in $ \hat{\mathcal
{L}}$. Let $\hat{K}$ be the unitary operator in $\hat{\mathcal{
L}}$, defined as
%
%e22 #&#
%
\begin{equation}
\label{coordwise} \hat{K} = I_{\hat{\mathcal{L}} \perp} + \sum
_{i=0}^{2\kappa+1}
a_i \langle b_i,\cdot\rangle_F ,
\end{equation}
where $I_{\hat{\mathcal{L}} \perp}$ is projector on the orthogonal
complement of $\hat{\mathcal{ L}}$ to $L_2(F)$. For convenience, let
us single out three short statements as a lemma.
%
%============================ Lemma == parametric == Unitarity
%==========
%
%le4 #&#

\begin{lemma}\label{lem:lemsec4}
\textup{(i)} The operator $\hat{K}$ is unitary on $\hat{\mathcal
{L}}$. It
maps basis $b$ into basis $a$ while it maps any function, orthogonal to
$\hat{\mathcal{L}}$ to itself:
\[
\hat{K}b=a, \quad\mbox{and}\quad \hat{K}\phi= \phi,\qquad
\mathit{if}\ \phi\perp\hat{\mathcal{L}}.
\]

\textup{(ii)} For a function $\phi$ consider its projection parallel to
functions $q_0,\dots, q_\kappa$,
\[
\phi- \sum_{i=0}^\kappa q_i
\langle q_i,\phi\rangle_F .
\]
Then
\[
\hat{v}_F(\phi) = \hat{v}_F\Biggl(\phi- \sum
_{i=0}^\kappa q_i\langle q_i,
\phi\rangle_F\Biggr) = w_F\Biggl( \phi- \sum
_{i=0}^\kappa q_i\langle q_i,
\phi\rangle_F\Biggr) .
\]
\end{lemma}

In other words, according to (ii), the processes $\hat{v}_F$
and $w_F$ coincide on the subspace of functions orthogonal to
$q_0,\dots, q_\kappa$. Both (i) and (ii) can be
easily checked. For example, the last equality follows from the
definition of $\hat{v}_F$ in Section~\ref{sec2}.
%
%============================Thm, parametric == the
%main==================
%
%th7 #&#

\begin{them}\label{thm:main2}
If $\hat{v}_F$ is $q$-projected $F$-Brownian motion and $G$ is
absolutely continuous with respect to $F$, then
\[
\hat{v}_G(\psi) = \hat{v}_F\bigl(\hat{K} (l\psi)\bigr)
\]
or, more explicitly,
%
%e23 #&#
%
\begin{equation}
\label{main2} \hat{v}_G(\psi) = \hat{v}_F(l\psi) -
\sum_{i=\kappa+1}^{2\kappa
+1} \hat{v}_F(a_i)
\langle l\psi, a_i-b_i\rangle_F
\end{equation}
is $r$-projected $G$-Brownian motion. If $G$ and $F$ are equivalent,
then this transformation is one-to-one.
\end{them}

From the point of view of this theorem, testing of various parametric
families with square integrable score functions of the same dimension
and equivalent $F_\theta$ and $G_{\theta'}$, is not a multitude of
various unconnected testing problems; since these testing problems can
be mapped into one another they can be glued in equivalence classes.
One representative from each class is, therefore, sufficient to use and
this makes the testing asymptotically distribution-free.
%
% The Theorem has surprising number of corollaries which look unusual
%and important in their own right. One hopes they are also useful in
%the goodness of fit theory.
%
\begin{pf*}{Proof of Theorem~\ref{thm:main2}} First, we prove that
$\hat{v}_G(\psi)$ is
$r$-projected $G$-Brownian motion in $\psi$, and then we show that
explicit expression of the right-hand side is that given in \eqref
{main2}. Consider
\begin{eqnarray*}
\hat{K} \phi& =&\phi- \sum_{i=0}^{2\kappa+1}
b_i \langle\phi, b_i \rangle_F + \hat{K}
\sum_{i=0}^{2\kappa+1} b_i \langle\phi,
b_i \rangle_F
\nonumber
\\
&=&\phi- \sum_{i=0}^{2\kappa+1} b_i
\langle\phi, b_i \rangle_F + \sum
_{i=0}^{2\kappa+1} a_i \langle\phi,
b_i \rangle_F .
\end{eqnarray*}
The second equality here uses part (i) of the lemma. The last display
in part (ii) shows that we need to consider projection of the
latter function parallel to $a_0,\dots, a_\kappa$. In taking this
projection, the sum $\sum_{i=0}^{\kappa} a_i \langle\phi, b_i
\rangle_F$ will be annihilated, so that the projection is
%
%e24 #&#
%
\begin{equation}
\label{project} \phi- \sum_{i=0}^{2\kappa+1}
b_i \langle\phi, b_i \rangle_F + \sum
_{i=\kappa+1}^{2\kappa+1} a_i \langle\phi,
b_i \rangle_F .
\end{equation}
Therefore, again using (ii),
\[
\hat{v}_F(\hat{K} \phi) = w_F\Biggl(\phi- \sum
_{i=0}^{2\kappa+1} b_i \langle\phi,
b_i \rangle_F + \sum_{i=\kappa+1}^{2\kappa+1}
a_i \langle\phi, b_i \rangle_F \Biggr).
\]
The first difference in \eqref{project} is orthogonal to the second
sum. Therefore,
\begin{eqnarray*}
E\hat{v}_F^2(\hat{K} \phi) &= &\langle\phi,\phi
\rangle_F - \sum_{i=0}^{2\kappa+1}
\langle\phi, b_i \rangle_F^2 + \sum
_{i=\kappa+1}^{2\kappa+1} \langle\phi, b_i
\rangle_F ^2
\\
& =& \langle\phi,\phi\rangle_F - \sum_{i=0}^{\kappa}
\langle\phi, b_i \rangle_F^2 .
\end{eqnarray*}
For $\phi=l\psi$, the latter expression is equal to
\[
\langle\psi,\psi\rangle_G - \sum_{i=0}^{\kappa}
\langle\psi, r_i \rangle_G^2 ,
\]
which is the variance of $\hat{v}_G(\psi)$. To arrive now at the
explicit form \eqref{main2} of $\hat{v}_G(\psi)$, rewrite \eqref
{project} as
\[
\phi- \sum_{i=0}^{2\kappa+1} a_i
\langle\phi, a_i \rangle_F + \sum
_{i=\kappa+1}^{2\kappa+1} a_i \langle\phi,
b_i \rangle_F
\]
and use orthogonality of $v_F$ to $a_0,\dots, a_\kappa$.
\end{pf*}

Weak convergence result, which follows from our theorem, is easy to
formulate in function-parametric as well as set-parametric versions,
but it is somewhat more convenient for application to consider, again,
the point-parametric version of the parametric empirical process, where
the family of functions $\psi$ is chosen as a family of indicator
functions, $\psi(y)=\mathbf{1}_{(\infty, x]}(y)$, indexed by~$x$.
Then transformation in \eqref{main2} applied to $v_n(\cdot, \hat
{\theta}_n)$ leads to the process
%
%e25 #&#
%
\begin{eqnarray}
\label{param-finite-n} \tilde v_{n} (x,\hat{\theta}_n) &= &\int
_{y\leq x} l(y) v_n(\mathrm{d}y,\hat{\theta}_n)\nonumber\\[-8pt]\\[-8pt]
&&{} -
\sum_{i=\kappa+1}^{2\kappa+1} \int_{y\in\Rd}
a_i(y) v_n(\mathrm{d}y,\hat{\theta}_n) \int_{y\leq x} \bigl(a_i(y)
- b_i(y)
\bigr) \,\mathrm{d}F(y) .
\nonumber
\end{eqnarray}
Weak convergence of the process $v_n(\cdot,\hat{\theta}_n)$ was
considered in a very large number of publications; among the first we
know of are \cite{KKW} and much later, but still long ago, \cite
{Dur73}. Certain (incomplete) review is given in \cite{ShW},
Chapter~3.5; convergence of $v_n(\phi,\hat{\theta}_n)$ on countably many
square integrable functions was studied in \cite{Khm79}. Based on
this, we take the weak convergence of the first integral in \eqref
{param-finite-n} as a process in $x$ as given, as well as convergence
of integrals from $a_i$ with respect to $v_n(\cdot,\hat{\theta}_n)$.
Their joint weak convergence is obvious and this leads to the statement
%
%e26 #&#
%
\begin{equation}
\label{empirical wc} \tilde v_n (\cdot,\hat{\theta}_n)
\stackrel{\mathcal{D}(F_\theta)} {\longrightarrow} \hat{v}_{G}
(\cdot) .
\end{equation}
If $F_\theta$ in our parametric family have rectangular support in
$\Rd$ then, as we already mentioned, the product distribution $G$
exists and we can proceed as in Corollary~\ref{cor:finite n}. However,
one point here needs some remark. The most natural choice of $G$ will
be a product of the marginal distributions of $F_\theta$ and,
therefore, $G=G_\theta$ will depend on $\theta$ as well. All
functions, $l=l_\theta, a_i=a_{i\theta}, b_i=b_{i\theta}$, which
participate in the transformation, will also depend on $\theta$. The
latter is true even if one chooses one common $G$ for all $F_\theta$,
simply because $F_\theta$ and, therefore, $\beta_{F_\theta}$ as
well, depend on~$\theta$. Hence, in \eqref{param-finite-n} the
functions $l_{\hat{\theta}_n}, a_{i \hat{\theta}_n}, b_{i \hat
{\theta}_n}$ will have to be used. This, however, creates only a minor
problem: in simple continuity assumptions on $l_\theta$ and $\beta
_{F_\theta}$ in $\theta$, similar, for example, to bracketing
assumptions in \cite{vdVW}, one can see that the difference between
transformation produced by $l_{\hat{\theta}_n}, a_{i \hat{\theta
}_n}, b_{i \hat{\theta}_n}$ and $l_\theta, a_{i\theta}, b_{i\theta
}$ is asymptotically small and, therefore, \eqref{empirical wc} is
still true.

More interesting and specific to this paper is the problem of practical
implementation and convenience of transformation \eqref{main2}.

%==========================Uniqueness==3.4===============
%s3.4 #&#
\subsection{Uniqueness of \texorpdfstring{$\hat{K}$}{hatK} and practical calculations of
\texorpdfstring{$v_n(\hat{K} \psi)$}{vn(hatKpsi)}}\label{sec3.4}

Start by noting that the operator $\hat{K}$ is an extension of the
operator $K$ of \eqref{unitary} to the parametric case. Moreover, the
former can be expressed by the latter. To show this, assume first
$\kappa=0$ and denote
\[
K_{g,h}=I - \frac{2}{\|h-g\|_F^2} (h-g)\langle h-g,\cdot
\rangle_F
\]
a unitary operator on $L_2(F)$ with the same properties as in
Lemma~\ref{lem:Koperat}, only with $l$ and $q_0$ replaced by general
$h$ and
$g$, respectively. Thus, $K_{q_0,l}=K$ of \eqref{unitary}. Now assume
$\kappa=1$. Recall that $K_{q_0,l}$ maps function $l$ to function
$q_0$ and maps any function, orthogonal to $l$ and $q_0$, to itself.
Consider the image of the function $\mathit{lr}_1$,
\[
K_{q_0,l}\mathit{lr}_1= \tilde{\mathit{lr}}_1.
\]
Since $l$ and $\mathit{lr}_1$ are orthogonal by construction, then so
are their
images $q_0$ and $\tilde{\mathit{lr}}_1$. Now consider operator
$K_{q_1, \tilde
{\mathit{lr}}_1}$. The product
%
%e27 #&#
%
\begin{equation}
\label{iterat} \hat{K} = K_{q_1, \tilde{\mathit{lr}}_1} K_{q_0,l}
\end{equation}
is another form of the operator $\hat{K}$. Indeed, as a product of
unitary operators, $\hat{K}$ is a unitary operator, and it maps any
functions $\phi_\perp$, which are orthogonal to $l, \mathit{lr}_1,
q_0, q_1$,
to itself, while it maps $l$ into $q_0$ and $\mathit{lr}_1$ into $q_1$:
\[
K_{q_1, \tilde{\mathit{lr}}_1} K_{q_0,l} l = K_{q_1, \tilde{\mathit
{lr}}_1} q_0 =
q_0 ,
\]
and
\[
K_{q_1, \tilde{\mathit{lr}}_1} K_{q_0,l} \mathit{lr}_1= K_{q_1,
\tilde{\mathit{lr}}_1} \tilde
{\mathit{lr}}_1 = q_1 .
\]
Since $b_2$ and $b_3$ are orthonormal and orthogonal to $l$ and
$\mathit{lr}_1$,
it follws that $\hat{K} b_2$ and $\hat{K} b_3$ also will be
orthonormal and orthogonal to $q_0$ and $q_1$, which is what is
required from $a_2$ and $a_3$.

This procedure can be iterated in $\kappa=2, 3$, and so forth. Hence,
it follows that transformation \eqref{main2} can be carried out as a
sequence of $\kappa+1$ just $one$-dimensional transformations. This
was tried recently in \cite{Thuong} with applications to testing
independence in contingency tables, and demonstrated that the coding is
simple and the calculations quick. In one of numerical examples, the
author considered $5\times6$ tables with, therefore, $\kappa=9$
marginal probabilities to estimate.

At the same time, comparison of the representation \eqref{iterat} with
the coordinate form used in Lemma~\ref{lem:lemsec4} rises the question
of uniqueness of $\hat{K}$, which is good to clarify.

To this end, consider the orthogonal decomposition of $L_2(F)$, which
uses the basis $b$:
\[
\hat{\mathcal{ L}}_\perp+ \hat{\mathcal{L}} = \hat{ \mathcal
{L}}_\perp+\mathcal{L} _{1b}+\mathcal{L}_{2b},
\]
where the subspace $\mathcal{L}_{1b}=\mathcal{L} (b_0,\dots,b_\kappa
)=\mathcal{L}(l,\ldots, \mathit{lr}_\kappa)$
is generated by the functions $b_0,\ldots,b_\kappa$, and $\mathcal
{L}_{2b}=\mathcal{L} (b_{\kappa+1},\ldots, b_{2\kappa+1})$ is
generated by the remaining part of the basis $b$, and $\hat{\mathcal
{L}}_\perp$ is the orthogonal complement of their sum to $L_2(F)$.
Similarly, consider orthogonal decomposition which uses the basis $a$:
\[
\hat{\mathcal{L}}_\perp+ \hat{\mathcal{L}} = \hat{\mathcal{
L}}_\perp+\mathcal{L} _{1a}+\mathcal{L}_{2a}.
\]
Then what the operator $\hat{K}$, defined in \eqref{coordwise}, does
is the following:
it maps unitarily subspace $\mathcal{L}_{ib}$ onto $\mathcal{L}_{ia},
i=1,2$, while leaves $\mathcal{L}_\perp$ unperturbed. However, let
$T_b$ be a unitary operator, which can be decomposed into direct sum
$T_b=T_\perp+T_{1b} + T_{2b}$ of unitary operators, of which $\hat
{\mathcal{L}}_\perp$, $\mathcal{L}_{1b}$ and $\mathcal{L}_{2b}$ are
invariant subspaces, respectively. Then, for any such operator, the
process
\[
v_{G,T_b}(\psi)=v_F(\hat{K} T_b l\psi)
\]
is also a $G$-Brownian motion. Moreover, if $T_a$ is a similar unitary
operator with invariant subspaces $\hat{\mathcal{L}}_\perp$,
$\mathcal{L}_{1a}$ and $\mathcal{L}_{2a}$, then
\[
v_{G,T_a,T_b}(\psi)=v_F(T_a \hat{K}
T_b l\psi)
\]
is again a projected $G$-Brownian motion. This makes nonuniqueness of
\eqref{coordwise} an obvious and, basically, trivial fact.

However, in practical problems we will not be in need to use $T_b$ and
$T_a$ in so much generality. Indeed, there does not seem to be a reason
to ``rotate'' $\phi_\perp$ and therefore we can agree to choose
$T_\perp$ as the identity operator on $\mathcal{L}_\perp$. Given
``target'' score functions, that is, given $l,\dots, \mathit
{lr}_\kappa$, and
the score functions $q_0, \dots, q_\kappa$ of the hypothetical
parametric family, it does not seem useful to ``rotate'' any of them and
one can agree to the rule that each $\mathit{lr}_i$ is mapped onto
$q_i$ for all
$i=0,\dots,\kappa$. This will uniquely define the image of $\phi
_{1b}$ as $\sum_{i=0}^\kappa a_i \langle\phi_{1b}, b_i\rangle$.
Moreover, for each $\phi$, the decomposition of $\phi=\phi_\perp
+\phi_{1b}+\phi_{2b}$ into its parts in the corresponding subspaces
is unique, and, in particular, $\phi_{2b}$ does not depend on the
choice of $b_i, i=\kappa+1,\dots, 2\kappa+1$, although the choice of
these latter functions is not unique.

More specifically, with the matrix
\[
C=\bigl\| \langle q_i, \mathit{lr}_j \rangle\bigr\|, \qquad i,j=0,\dots,
\kappa
\]
the coordinate functions of the vector
\[
(q_0,\dots, q_\kappa)^T - C (l,
\dots,\mathit{lr}_\kappa)^T
\]
are orthogonal to coordinates of $(l,\dots,\mathit{lr}_\kappa)^T$ and,
therefore, the vector $(b_{\kappa+1},\dots, b_{2\kappa+1})^T$ has
to be a linear transformation of the latter:
\[
(b_{\kappa+1},\dots, b_{2\kappa+1})^T = H
\bigl[(q_0,\dots, q_\kappa)^T - C (l,
\dots,\mathit{lr}_\kappa)^T\bigr] .
\]
This linear transformation $H$ renders the coordinates of $(b_{\kappa
+1},\dots, b_{2\kappa+1})^T$ mutually orthogonal and normalized.
However, the $H$ is not defined uniquely. Therefore, although with our
agreement, the vector
\[
\hat{K} \bigl[(q_0,\dots, q_\kappa)^T - C (l,
\dots,\mathit{lr}_\kappa)^T\bigr]=(l,\dots,\mathit{lr}_\kappa)^T
- C(q_0,\dots, q_\kappa)^T
\]
remains the same for any choice of operator $\hat{K}$ with properties
as in Lemma~\ref{lem:lemsec4}, nonuniqueness of $H$ makes the
multiple choice of $a_{\kappa+1},\dots, a_{2\kappa+1}$ possible.

Apart from simplicity in numerical calculations, the advantage of the
representation \eqref{iterat} is that it offers a unique
``canonical'' form of transformation. Then there is no need to be
interested in the form of $a_{\kappa+1},\dots, a_{2\kappa+1}$, as
they do not enter in our transformation $\hat{v}_G=\hat{v}_F(\hat{K}
l\psi)$ explicitly.

%==========================Numerical Illustrations===============

%s4 #&#
\section{Some numerical illustrations}\label{sec4}

Let $\tilde u_n$ denote the process obtained as transformation \eqref
{unif} applied to empirical process $v_{nF}$:
\[
\tilde u_n (x) = \int_{y\leq x} \frac{1}{\sqrt{f(y)}}
v_{nF}(\mathrm{d}y) - \frac{\int_{y\leq x}(1- \sqrt{f(y)} ) \,
\mathrm{d}y}{1- \int_{[0,1]^d} \sqrt
{f(y)} \,\mathrm{d}y}\int_{[0,1]^d}
\frac{1}{\sqrt{f(y)}} v_{nF}(\mathrm{d}y).
\]
The choice of $d=1$ suggested itself by the fact that the limit
distributions of statistics below are known and, therefore, one can
easily judge how quick is the convergence.
%f1 #&#
%
\begin{figure}

\includegraphics{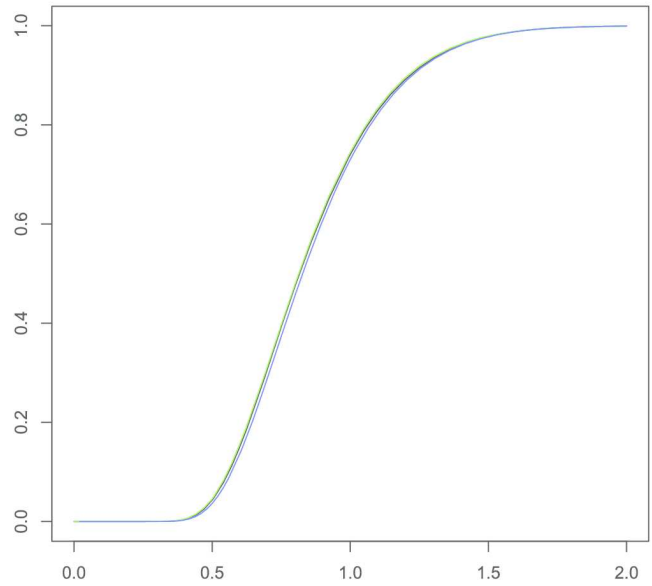}

\caption{Distribution functions of K-S statistics $D(\tilde u_n)$ for
the beta distributions, with bell-shaped and J-shaped densities,
described in the text. We used 10\,000 simulations of samples of size
$n=200$. The third is the graph of Kolmogorov distribution, which is
their limit in $n$.}
\label{fig:contin}
\end{figure}

In Figure~\ref{fig:contin}, two distribution functions of the statistic
\[
D(\tilde u_n) = \sup_{0<x<1} \bigl|\tilde
u_n(x)\bigr|
\]
are shown, for sample size $n=200$. It is not easy to distinguish them,
although the statistics are based on samples from quite different beta
distributions: with a bell-shaped (parameters 3 and 3) and J-shaped
(parameters 0.8 and 1.5) beta densities, respectively. The third graph
is that of the Kolmogorov distribution function, which is the limiting
distribution of the $D(u)=\sup_{0<x<1} |u(x)|$. If $\tilde u_n$ were a
sort of an empirical process, like, say $v_{n\tilde F}$ with some
$\tilde F$, the distribution of its supremum will again be that of
$D(u)$ and some doubts would remain whether $\tilde u_n$ behaves as a
uniform empirical process or an empirical process based on some other
distribution. However, our $\tilde u_n$ is not an empirical process at
all -- it is a difference between some weighted version of an empirical
process and some deterministic function times a linear functional from
the former.
%
%f2 #&#
%
\begin{figure}

\includegraphics{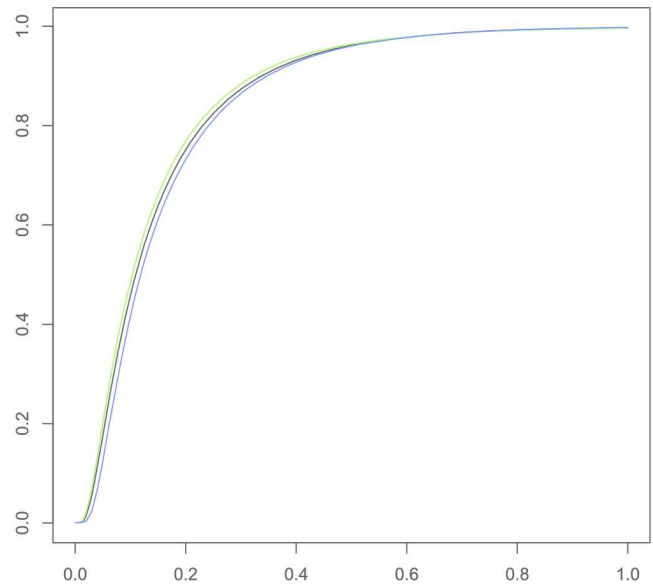}

\caption{Distribution functions of $\omega^2$-statistics for the same
underlying beta distributions as above. We used 10\,000 simulations of
samples of reduced size $n=50$. The lowest is the graph of $\omega^2$
distribution, which is their limit in $n$.}\vspace*{-10pt}
\label{fig:omega}
\end{figure}

Now, on Figure \ref{fig:omega}, we show distribution functions of the omega-square statistic
\[
\Omega^2 = \int_0^1 \bigl[\tilde
u_n(t)\bigr]^2 \,\mathrm{d}t .
\]
These distribution functions cannot converge to the omega-square
distribution unless $\tilde u_n$ indeed behaves as the uniform
empirical process. But it seems that they do. Although the differences
are now visible, note that the integral was calculated merely as a
Darboux sum with not too fine step, and that the sample size was only
$n=50$.

It is interesting to have some indication of how quickly the processes
of Theorem~\ref{cor:nine}
converge to Brownian motion. The point of particular interest was
whether division by $\sqrt f$, as in \eqref{bm21}, spoils the
convergence, and if so, by how much. For this comparison we used still
another version of $b_n$, which one obtains by integrating $\sqrt{f}$
with respect to the process \eqref{bm2}. For one-dimensional time, it
leads to
%
%e28 #&#
%
\begin{equation}
\label{bm3} \int_\Delta^x\sqrt{f(y)}
b_n(\mathrm{d}y) = v_{nF}(x)-v_{nF}(\Delta) +
v_{nF}(\Delta) \frac{F(x) - F(\Delta)}{\sqrt{ F(\Delta)} -
F(\Delta)} ,
\end{equation}
which certainly converges as quickly as empirical process $v_{nF}$.

\begin{figure}

\includegraphics{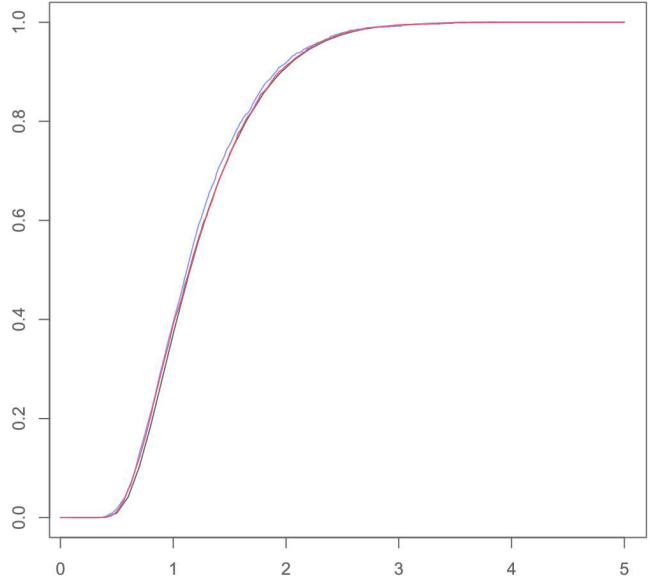}

\caption{Distribution functions of statistics $D(b_n)$ from the
process \protect\eqref{bm21} for the same underlying beta
distributions, as
above. Again, 10\,000 simulations of samples of size $n=200$. The third
is the graph of the distribution of $\sup_{0<x<1} |b(x) |$, which is
the limit distribution for the first two.}
\label{fig:normedmotion}
\end{figure}

Figure~\ref{fig:normedmotion} shows the graphs of distribution
function of K-S statistic from the process \eqref{bm21},
\[
D(b_n) =\sup_{\Delta<x<1} \bigl|b_n(x) \bigr|/{\sqrt{1-
\Delta}},
\]
obtained for two different beta distributions (described above) along
with the distribution function of supremum of a standard Brownian
motion. We see that the discrepancy between pre-limiting distribution,
for $n=200$, and the limit exists, but is very small, especially if we
consider convergence of quantiles.

%f3 #&#
%

%f4 #&#
%
\begin{figure}

\includegraphics{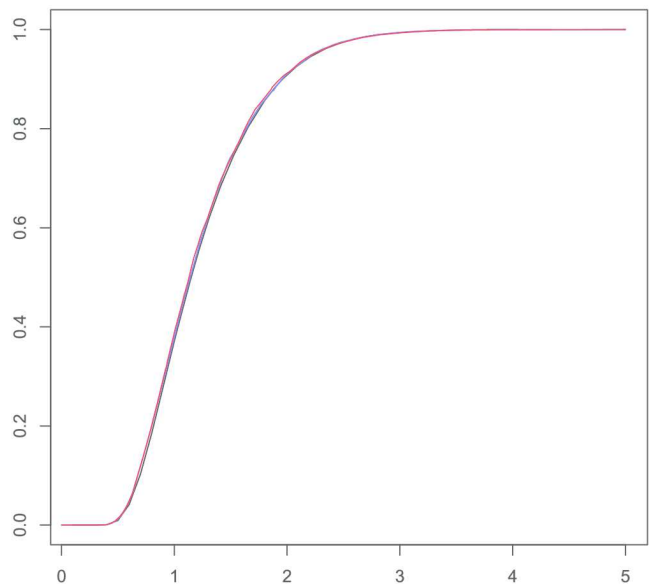}

\caption{Distribution functions of K-S statistic from the process
\protect\eqref{bm3} for the same underlying beta distributions as
above, and
with 10\,000 simulations of samples of size $n=200$.}
\label
{fig:motion}
\end{figure}

The last Figure~\ref{fig:motion} shows distribution functions of K-S
statistic from the process \eqref{bm3} normalized by $\sqrt
{1-F(\Delta)}$ for samples from the same underlying distributions as
in Figure~\ref{fig:normedmotion}. With respect to the previous figure,
there is some improvement, but not by much.

%\begin{appendix}
%\section{}
%\end{appendix}

\section*{Acknowledgements}
I want to thank the colleagues with whom I had possibility to discuss
the topic of this paper, most notably Professor P. Greenwood, as well
as Professor N. Henze and Professor H.L. Koul. I also thank the
associate editor and the referee for the high quality of their reading
and comments.

%\begin{supplement}%[id=suppA]
%\slink[doi]{10.3150/00-BEJXXXXSUPP} %[doi,text={...}] - jei reikia
%suskaldyti doi
%\sfilename{BEJ000\_supp.pdf}
%\end{supplement}

%
% imsref loaded by jurgita.kaciuliene, 2014-08-20 11:10:38

\printhistory
\end{document}